\documentclass[12pt]{article}
\usepackage{fullpage,amsmath,amssymb,amsthm,hyperref}
\everymath{\displaystyle}
\allowdisplaybreaks[3]

\newcommand{\Z}{\mathbb{Z}}

\renewcommand{\P}{p^{\textup{pk}}}
\newcommand{\PP}{P^{\textup{pk}}}
\newcommand{\pk}{\textup{pk}}

\begin{document}
\title{Period mimicry: \\ A note on the $(-1)$-evaluation of the peak polynomials}
\author{Justin M. Troyka \\ York University, Toronto}
\maketitle

Let $\P_n(t)$ be the peak polynomial for permutations of size $n$, i.e.\ $\P_n(t) = \sum_{\pi \in S_n} t^{\pk(\pi)}$ (where $\pk(\pi)$ is the number of peaks of $\pi$). Define
\[ \PP(z,t) = \sum_{n\ge0} \P_n(t) \frac{z^n}{n!} \quad \text{and} \quad F(z) = \PP(z,-1), \]
and let $f_n$ denote the coefficient of $z^n/n!$ in $F(z)$, namely $f_n = \P_n(-1)$. The sequence $f_n$ is the object of study in this note.

It is well known that
\[ \PP(z,t) = \frac{\sqrt{1-t}}{\sqrt{1-t} - \tanh(z\sqrt{1-t})}, \]
whence
\begin{equation} \label{eq:formula} F(z) = \frac{\sqrt{2}}{\sqrt{2} - \tanh(z\sqrt{2})}. \end{equation}
Thus $f_n$ is the coefficient of $z^n/n!$ in $\frac{\sqrt{2}}{\sqrt{2} - \tanh(z\sqrt{2})}$, which is given by \href{http://oeis.org/A006673}{\texttt{A006673}} in the OEIS \cite{OEIS}. Tirrell and Zhuang \cite{TZ} study this sequence from the point of view of bijective combinatorics. The purpose of this note is to study it from the point of view of analytic combinatorics, from which we discover a curious phenomenon.

The first few values of $f_n$ (starting at $n=0$) are:
\[ 1,\, 1,\, 2,\, 2,\, -8,\, -56,\, -112,\, 848,\, 9088,\, 25216,\, -310528,\, -4334848,\, -14701568,\,\ldots \]
One might conjecture that, for $n \ge 1$, $f_n$ is positive if $(\text{$n$ mod $6$}) \in \{1,2,3\}$ and negative if $(\text{$n$ mod $6$}) \in \{4,5,0\}$. However, the following analysis shows that the sequence's oscillation is irregular, and indeed if we compute enough terms of the power series we can see that
\[ f_{42} = 356\,077\,960\,394\,850\,110\,410\,690\,594\,606\,123\,271\,850\,033\,152 \approx 3.56 \times 10^{44} \]
breaks the apparent pattern of signs (and is the first term to do so).

The asymptotic behavior of $f_n$ is determined by the location and nature of the dominant singularities of $F(z)$, considered as a function of a complex variable. (The \emph{dominant singularities} are the singularities closest to $0$; they lie on the boundary of the disk of convergence.) What follows is an asymptotic analysis of $f_n$. To skip right to the answer, go to Section \ref{sec:answer}.

\section{Asymptotic analysis}

The singularities of $F(z)$ occur where the bottom of the fraction in \eqref{eq:formula} equals $0$, i.e.\ $\sqrt{2} - \tanh(z\sqrt{2}) = 0$. We solve for $z$ and obtain
\[ z = \frac{\log(3+2\sqrt{2}) + (2k+1)\pi i}{2\sqrt{2}} \quad (k \in \Z), \]
of which there are two that are closest to $0$:
\[ z = \frac{\log(3+2\sqrt{2}) \pm \pi i}{2\sqrt{2}}. \]
The fact that there are two dominant singularities means we should expect $f_n$ to exhibit some kind of oscillatory behavior.

Define $z_0 = \frac{\log(3+2\sqrt{2}) + \pi i}{2\sqrt{2}}$, so the two dominant singularities of $F(z)$ are $z_0$ and $\overline{z_0}$. The asymptotic behavior of $f_n$ comes from the behavior of $F(z)$ near $z_0$ and $\overline{z_0}$. We can decompose $F(z)$ as
\begin{equation} \label{eq:decomp} F(z) = \frac{\sqrt{2}}{\sqrt{2} - \tanh(z\sqrt{2})} = \frac{-1}{z_0} A(z) + \frac{-1}{\left(\overline{z_0}\right)} B(z),  \end{equation}
where
\[ A(z) = \frac{2(z-z_0)(z-\overline{z_0})}{\pi i (\sqrt{2} - \tanh(z\sqrt{2}))} \cdot \frac{1}{1-z/z_0} \quad \text{and} \quad B(z) = -\frac{2(z-z_0)(z-\overline{z_0})}{\pi i  (\sqrt{2} - \tanh(z\sqrt{2}))} \cdot \frac{1}{1-z/\overline{z_0}}. \]
First consider $A(z)$. We have $A(z) \sim \frac{1}{1-z/z_0}$ as $z \to z_0$, in the sense that $\lim_{z\to z_0} \frac{A(z)}{1/(1-z/z_0)}  = 1$. Therefore, since $A(z)$ is analytic on the domain $\{\text{$|z/z_0| < 1 + \varepsilon$ and $z \not= z_0$}\}$ (for some $\varepsilon > 0$), we can apply \cite[Cor.\ VI.1, p.\ 392]{FS} and obtain
\begin{equation} \text{$[z^n] A(z) \sim {\left(1/z_0\right)}^n$ as $n \to \infty$}. \label{eq:A} \end{equation}
Furthermore, since $A(z) - \frac{1}{1-z/z_0}$ is analytic on the domain $\{|z/z_0| < 1 + \varepsilon\}$, the asymptotic expansion in \eqref{eq:A} has an exponentially decaying error term:
\[ \text{$[z^n] A(z) = {\left(1/z_0\right)}^n + O(\alpha^n)$ as $n \to \infty$} \]
for some $\alpha$ satisfying $0 < \alpha < 1/|z_0|$.

The same argument with $B(z)$ shows that
\[ \text{$[z^n] B(z) = {\left(1/\overline{z_0}\right)}^n + O(\alpha^n)$ as $n \to \infty$}. \]
We now take these asymptotic formulas for $[z^n]A(z)$ and $[z^n]B(z)$ and put them into \eqref{eq:decomp}:
\begin{align*}
f_n/n! = [z^n]F(z) &= \frac{-1}{z_0} \,[z^n]A(z) + \frac{-1}{\left(\overline{z_0}\right)} \,[z^n]B(z) \\
&= -\left[{\left(z_0\right)}^{-(n+1)} + {\left(\overline{z_0}\right)}^{-(n+1)}\right]  + O(\alpha^n).
\intertext{Furthermore, letting $\theta$ denote the angle of $z_0$ from the positive real axis, we obtain $z_0 = \left|z_0\right| e^{i\theta}$ and $\overline{z_0} = \left|z_0\right| e^{-i\theta}$, and so}
f_n/n! &= -\left[\left|z_0\right|^{n+1} e^{-i(n+1)\theta} + \left|z_0\right|^{n+1} e^{i(n+1)\theta} \right] + O(\alpha^n) \\
&= -2\,{\left|z_0\right|}^{-(n+1)} \cos((n+1)\theta) + O(\alpha^n).
\end{align*}

\section{Conclusion} \label{sec:answer}
Those computations show that
\begin{gather*}
f_n = -2\,\rho^{-(n+1)} \cos((n+1)\theta) \cdot n! + O(\alpha^n n!), \\
\text{where $\rho = |z_0| = \frac{\sqrt{\left[\log(3+2\sqrt{2})\right]^2 + \pi^2}}{2\sqrt{2}} \approx 1.274$} \\
\text{and $\theta = \arg(z_0) = \arctan\left(\frac{\pi}{\log(3+2\sqrt{2})}\right) \approx 1.012 \cdot \pi/3$} \\
\text{and $\alpha$ is some number in $(0,\rho^{-1})$}.
\end{gather*}

We remark that this asymptotic formula makes it immediately clear why $f_n$ at first appears to exhibit periodic behavior. If $\theta$ were exactly equal to $\pi/3$, then the factor of $\cos((n+1)\theta)$ would result in precisely the $6$-periodic sequence of signs that $f_n$ seems to have for $n < 42$. Since $\theta$ is very close to $\pi/3$ but not equal to it, the sequence gradually drifts away from this $6$-periodic pattern. This is a rare sighting of period mimicry in the wild.

Does our asymptotic formula imply that $f_n$ has the same sign as $-\cos((n+1)\theta)$ for large enough $n$? Not quite: $\cos((n+1)\theta)$ may sometimes get very close to $0$, and if $\cos((n+1)\theta)$ is small compared to $(\alpha \rho)^n$ then the asymptotic expansion may be dominated by the error term $O(\alpha^n n!)$. However, because $(\alpha \rho)^n$ is decaying exponentially, this almost never occurs (in the sense that, as $n \to \infty$, the proportion of $k \le n$ for which the error term dominates goes to $0$). Thus $f_n$ and $\cos((n+1)\theta)$ ``usually'' have the same sign, and indeed it may well be true that they have the same sign for all $n \ge 1$.

The fact that $\theta$ is not a rational multiple of $\pi$ means that there is no $k$ for which the sequence's fluctuations are $k$-periodic (see \cite[``Nonperiodic fluctuations'', pp.\ 264--266]{FS}). Since the sign of $f_n$ is unpredictable, it seems there is little hope of finding a combinatorial interpretation of $f_n$ in which permutations with opposite signs cancel each other out.

\subsection*{Acknowledgements}
I would like to thank Jordan Tirrell for introducing me to this problem and for helpful discussions. I would like to thank Yan Zhuang for asking me whether the asymptotic formula implies that $f_n$ has the same sign as $-\cos((n+1)\theta)$, which prompted me to more carefully carry out the analysis and to add the paragraph that addresses this question.

\end{document}